\documentclass[11pt]{amsart}
\usepackage{amssymb}
\usepackage{amscd}
\usepackage{curves}
\usepackage{epsfig}
\usepackage[pass]{geometry}
\usepackage{graphicx}
\usepackage{pst-grad}
\usepackage{pst-plot}
\usepackage{verbatim}
\usepackage{latexsym,eucal,amsfonts,amssymb,amsmath,graphicx}
\usepackage{caption}
\numberwithin{equation}{section}
\newtheorem{theorem}{Theorem}[section]

\newtheorem{prop}[theorem]{Proposition}

\newcommand\fnote[1]{\captionsetup{font=small}\caption*{#1}}
\def \bpf {\begin{proof}}
\def \epf {\end{proof}}
\def \beq {\begin{equation*}}
\def \eeq {\end{equation*}}
\def \bsp{\begin{split}}
\def \esp{\end{split}}
\def \beqq {\begin{equation}}
\def \eeqq {\end{equation}}

\def \mcd {{\mathcal D}}
\def \mce {{\mathcal E}}

\def \mbr {{\mathbb R}}

\def \diag{\textrm{Diag}}

\def \supp {\text{supp }}

\def\Id {\operatorname{Id}}

\def \eps {\epsilon}   
\def \la {\lambda}   
\def \La {\Lambda}

\def \lap {\Delta}
\def \p {\partial}
\def \ha {\frac{1}{2}}

\def \WF {\operatorname{WF}}
\def \singsupp {\operatorname{singsupp}}

\begin{document}
\title{Reducing streaking artifacts in quantitative susceptibility mapping}
\author{Benjamin Palacios	}
\address{Benjamin Palacios
\newline
\indent Department of Mathematics, University of Washington}
\email{bpalacio@uw.edu}

\author{Gunther Uhlmann} 
\address{Gunther Uhlmann
\newline
\indent Department of Mathematics, University of Washington,
\newline
\indent Institute for Advanced Study, the Hong Kong University of Science and Technology 
\newline
\indent and  Department of Mathematics, University of Helsinki}
\email{gunther@math.washington.edu}

\author{Yiran Wang}
\address{Yiran Wang
\newline
\indent Department of Mathematics, University of Washington,
\newline
\indent and Institute for Advanced Study, the Hong Kong University of Science and Technology}
\email{wangy257@math.washington.edu}

\begin{abstract} 
It is well-known that reconstruction algorithms in quantitative susceptibility mapping often contain streaking artifacts. In \cite{Seo}, the cause of the  artifacts is identified as propagation of singularities. In this work, we analyze such singularities carefully and propose some strategies to reduce the artifacts. 
\end{abstract}

\maketitle

\section{The inverse problem}
The goal of  quantitative susceptibility mapping (QSM) is to provide images of the magnetic susceptibility distribution $\chi$ inside the human body. We refer to \cite{SeWo, Wa} for the detailed background. Roughly speaking, when a tissue is put in a known magnetic field $\mathbf{H}$, it acquires a magnetic moment $\mathbf{M}$. The magnetic susceptibility $\chi$ of the tissue is defined by $\mathbf{M} = \chi \mathbf{H}$ where $\chi$ represents the electronic perturbation in the tissue.  QSM aims to reconstruct $\chi$ from the measured local field perturbations. Actually, the experiment is carried out in a magnetic resonance (MR) scanner. The tissue is put in a known magnetic field and we measure the field $\psi$ associated with the magnetization polarized by the main magnetic field of the MR scanner. Mathematically, this process can be described by   
\beqq\label{qsm1}
\psi(x) = \text{p.v.} \int_{\mbr^3} d(x -x') \chi(x') dx',
\eeqq
where $x = (x_1, x_2, x_3), x' = (x_1', x_2', x_3') \in \mbr^3$ and  $\text{p.v.}$ denotes the principal value of the singular integral with kernel 
\beq
d(x) = \frac{2x_3^2 - x_1^2 - x_2^2}{4\pi |x|^5}, \ \ x\in \mbr^3.
\eeq
The inverse problem is to find $\chi$ given $\psi$, which is a deconvolution problem.

It is convenient to work in the Fourier domain. Denote the Fourier transform of $f$ by $\hat f$. By taking the Fourier transform of \eqref{qsm1}, we obtain that
\beqq\label{qsm2}
\hat \psi(\xi) = D(\xi) \hat \chi(\xi), \ \ \xi = (\xi_1, \xi_2, \xi_3)\in \mbr^3,
\eeqq
where 
\beq
D(\xi) = \frac{1}{3} - \frac{\xi_3^2}{|\xi|^2}.
\eeq  
The inverse problem of QSM is equivalent to recover $\hat \chi$ from $\hat \psi$. It is clear that the problem is ill-posed due to the presence of  zeros of $D(\xi)$. In practice, since the data $\psi$ may contain errors, it is noticed that the reconstructed image of $\chi$ is often contaminated by streaking artifacts, see \cite{Seo} and the references therein. For example, see Figure \ref{fig:phantom_data_seo}. There have been several works on removing the artifacts, mostly based on regularization techniques, see \cite{SeWo, Wa} for a review and references. In \cite{Seo}, the authors performed the mathematical study of QSM. For $\psi$ in a suitable distribution space, the existence, uniqueness and reconstruction formulas are obtained in \cite[Theorem 2.2]{Seo}. Also, the cause of the streaking artifacts are identified as propagation of singularities for wave-type operators. In this note, we analyze the streaking artifacts in more detail. We will give a more precise description of the singularities under some assumptions and propose some strategies to reduce the streaking artifacts. The microlocal analysis of the singularities involves paired Lagrangian distributions that were introduced in \cite{MU} and further studied on \cite{DUV, FLU, GrU0, GrU, GrU1, GS, GU}. In particular, the method developed in this paper could be applied to reduce artifacts for restricted X-ray transforms \cite{GrU, FLU} like the X-ray transform with sources on a curve. 


\section{Streaking artifacts and propagation of singularities}

As discussed in \cite[Section 2.2]{Seo}, equation \eqref{qsm2} is equivalent to the wave equation on $\mbr^2$. Actually, multiplying \eqref{qsm2} by $|\xi|^2$ and taking the inverse Fourier transform, we obtain that 
\beqq\label{eqwave}
P(\p)\chi(x) = -\lap \psi(x),
\eeqq
where $P(\p)$ is the wave-type operator 
\beqq\label{op_P}
P(\p) = 
\frac{2}{3}\frac{\p^2}{\p x_3^2} - \frac{1}{3}(\frac{\p^2}{\p x_1^2} + \frac{\p^2}{\p x_2^2}),
\eeqq
and $\lap = \sum_{i = 1}^3\frac{\p^2}{\p x^2_i}$ is the Laplacian on $\mbr^3$. It is well-known that the wave operator $P(\p)$ has a fundamental solution i.e. there exits $Q\in \mcd'(\mbr^3\times \mbr^3)$, the space of distributions on $\mbr^3\times \mbr^3$, such that 
\beq
P(\p) Q = \delta,
\eeq
where $\delta$ is the delta distribution on $\mbr^3.$ In fact, one can write down an expression of $Q$ explicitly i.e.\ $Q(x, x') = g(x-x')$ where 
\beqq\label{fund_sol}
g(x) = 
\left\{\begin{array}{c}
\cfrac{3}{4\pi \sqrt{x_3^2 - 2(x_1^2 + x_2^2)}},  \text{ for } 2(x_1^2 + x_2^2) < x_3^2, \\ [1em]
0, \text{ otherwise}
\end{array}\right.
\eeqq
see (3.2) in \cite{Seo}. $Q$ defines a continuous linear map from $C_0^\infty(\mbr^3)$ to $\mcd'(\mbr^3)$ (by e.g.\  Schwartz kernel theorem). Hereafter we do not distinguish the notation for the operator and its Schwartz kernel. For $\psi\in C_0^\infty(\mbr^3)$, we can express the solution $\chi$ to \eqref{eqwave} simply as
\beqq\label{defQ}
\chi  = Q(-\lap \psi) = -g\ast (\lap \psi),
\eeqq
which is a convolution. The extension of $Q$ to distributions was done in \cite{Seo}. Here we give an exposition using microlocal techniques. Also, we give a detailed description of the singularities of $\chi$. 

It is easy to see that the singular support of $Q$ is the set $\{x\in \mbr^3: 2(x_1^2 + x_2^2) = x_3^2\}$. Recall that for any distribution $\phi$, the singular support of $\phi$, denoted by $\singsupp(\phi)$, is defined as the closure of the complement of the set where $\phi$ is smooth. To get a more precise description of the singularities, we use the notion of wave front sets (see for example \cite{Du, Ho1}) defined on the cotangent space $T^*\mbr^3$ which can be identified with the product space $\mbr^3\times \mbr^3$.  For $\phi\in \mcd'(\mbr^3)$, the wave front set $\WF(\phi)$  is a conic set in $T^*\mbr^3\backslash 0$ and by definition, $(x', \xi')\notin \WF(\phi)$ if there exists a conic neighborhood $\Gamma\subset T^*\mbr^3\backslash 0$ of $(x', \xi')$ such that for any $N>0$, there is $C_N>0$ such that
\beq
|\hat \phi(\xi)| \leq C_N |\xi|^{-N}, \ \ (x, \xi)\in \Gamma.
\eeq
It is a fact that $\pi(\WF(\phi)) = \singsupp(\phi)$ where $\pi: T^*\mbr^3\rightarrow \mbr^3$ denotes the natural projection $ \pi(x, \xi) = x$. 

The fundamental solution $Q$ belongs to a special class of distributions, namely the paired Lagrangian distributions. In particular, the wave front set of $Q$ consists of two intersecting Lagrangian submanifolds. The theory for such distributions is developed in \cite{DUV, MU, GU, GrU, GrU1}. It is now widely used in microlocal analysis and inverse problems, see for example \cite{FLU, GrU0, KLU, KLU1, LUW}.

We describe the two Lagrangians first. Let $\xi =  (\xi_1, \xi_2, \xi_3) \in T^*_x\mbr^3$ be the dual variables to $x =  (x_1, x_2, x_3) \in \mbr^3$. It is well-known that $T^*\mbr^3$ is a symplectic manifold with canonical two form given by 
\beq
\omega = d\xi\wedge dx = \sum_{i = 1}^3 d\xi_i\wedge dx_i.
\eeq
A submanifold $\La \subset T^*\mbr^3$ is called Lagrangian if $\dim \La = 3$ and the canonical two form $\omega$ vanishes on $\La$. A Lagrangian manifold $\La$ is conic if $(x, \la \xi)\in \La$ for any $(x, \xi)\in \La$ and $\la > 0$. Let $p(\xi)$ be the symbol of $P(\p)$ i.e.\ $p(\xi) = -\xi_3^2 + \frac{1}{3}|\xi|^2$. The characteristic set of $p$ is $\Sigma = \{(x, \xi)\in T^*\mbr^3: p(\xi) = 0\}.$ We denote the Hamilton vector field of $p$ by $H_p$. Explicitly, we have 
\beq
H_p = \sum_{i = 1}^3 (\frac{\p p}{\p \xi_i} \frac{\p }{\p x_i} - \frac{\p p}{\p x_i}\frac{\p }{\p \xi_i}) = \frac{2}{3}(\xi_1\frac{\p}{\p x_1} + \xi_2\frac{\p }{\p x_2})  -\frac{4}{3}\xi_3 \frac{\p }{\p x_3}.
\eeq
Notice that $H_p$ is tangent to $\Sigma$. The integral curves of $H_p$ in $\Sigma$ are called null bicharacteristics. It is a well-known fact that the projections of null bicharacteristics to the base manifold are geodesics. In our case, the base manifold is $\mbr^3$ with Lorentzian metric $g = -\frac{3}{2}dx_3^2 + 3 (dx_1^2 + dx_2^2)$ so the geodesics are straight lines. More precisely, for any $(x', \xi')\in \Sigma$, we denote the null bicharacteristics by $\gamma_{x', \xi'}(s), s\in \mbr$. Then we find that 
\beq
\gamma_{x', \xi'}(s) = (x' + s dp(\xi'), \xi'),
\eeq
where $dp(\xi) = (\frac{2}{3}\xi_1, \frac{2}{3}\xi_2, -\frac{4}{3}\xi_3)$. The projection of $\gamma_{x', \xi'}$ to $\mbr^3$ is just $x'+ s dp(\xi')$ which is a straight line. 

We also work on the product space $\mbr^3\times \mbr^3$ since $Q$ is a distribution defined there. We can regard the symbol $p$ as a function on $T^*\mbr^3\times T^*\mbr^3$ by lifting it from the left factor i.e.\ $p(\xi, \xi') = p(\xi)$. Similarly, we identify $\Sigma, H_p$ as objects on $T^*\mbr^3\times T^*\mbr^3$. We let 
\beq
\diag = \{(x, x')\in \mbr^3\times \mbr^3: x = x'\}
\eeq
be the diagonal of the product space and  
\beq
N^*\diag = \{(x, x'; \xi, \xi')\in T^*\mbr^3\times T^*\mbr^3: x = x', \xi' = -\xi, \xi\neq 0\}
\eeq 
be the conormal bundle of $\diag$. This is a conic Lagrangian submanifold in $T^*(\mbr^3\times \mbr^3)$ with canoic two form $\tilde w =  d\xi\wedge dx +  d\xi'\wedge dx'$, and gives the normal direction to $\diag$. Let $\La_p$ be the conic Lagrangian manifold in $T^*(\mbr^3\times \mbr^3)$ obtained from flowing out $N^*\diag\cap \Sigma$ under $H_p$. Actually, this can be written down explicitly as
\beqq\label{flowout}
\La_p = \{(x' + sdp(\xi'), \xi', x', -\xi')\in T^*\mbr^3\backslash 0 \times T^*\mbr^3\backslash 0: x'\in \mbr^3, s\in \mbr,\ \ p(\xi') = 0 \}.
\eeqq
It is proved in \cite{MU} that  the fundamental solution $Q\in I^{-\frac{3}{2}, -\ha}(N^*\diag, \La_p)$. We will explain this notion especially the meaning of the orders in the next section. For now, we just need the fact that the wave front set $\WF(Q)\subset N^*\diag\cup \La_p$. 

Let's recall the wave front relation for $E\in \mcd'(\mbr^3\times \mbr^3)$ which is
\beq
\WF'(E) = \{((x, \xi), (x', \xi'))\in T^*\mbr^3\times T^*\mbr^3\backslash 0: (x, x'; \xi, -\xi')\in \WF(E) \}.
\eeq
For two wave front relations $R_1, R_2 \subset T^*\mbr^3\times T^*\mbr^3$, the composition is defined as
\beq
R_1\circ R_2 = \{((x, \xi), (x'', \xi'')): \exists (x', \xi')\in T^*\mbr^3 \text{ s.t. } ((x, \xi), (x', \xi'))\in R_1, ((x', \xi'), (x'', \xi''))\in R_2\}.
\eeq
Now we have the following result on the solvability of QSM (compare with Theorem 3.3 of \cite{Seo}).
\begin{prop}\label{wf}
$Q$ can be extended to a sequentially continuous mapping: $\mce'(\mbr^3)\rightarrow \mcd'(\mbr^3)$.  For $\psi\in \mce'(\mbr^3)$, there exists a solution $\chi = Q(-\lap\psi)\in \mcd'(\mbr^3)$ to \eqref{eqwave} and
\beq
\WF(\chi)\subset \WF(\psi)\cup (\La_p'\circ \WF(\psi)).
\eeq
Moreover, for any $(x, \xi)\in (\WF(\chi)\backslash\WF(\psi)) \cap \Sigma$, let $\gamma_{x, \xi}$ be the bicharacteristics from $(x, \xi)$. Then $\gamma_{x, \xi} \subset \WF(\chi).$
\end{prop}
\bpf
The extension of $Q$  to $\mce'(\mbr^3)$  and the wave front relation are direct consequences of Corollary 1.3.8 of \cite{Du}. Here we used $\WF(\lap \psi) = \WF(\psi)$ because $\lap$ is an elliptic differential operator, see for example Corollary 8.3.2 of \cite{Ho1}. The last conclusion follows from the standard propagation of singularities result for wave operators (more generally operators with real principal part), see for example Theorem 8.3.3 of \cite{Ho3}. 
\epf

We remark that among the wave front set of $Q$, $N^*\diag$ does not move the singularities of $\psi$, however, $\La_p$ does. We know a priori that $\chi$ is compactly supported hence $\singsupp(\chi)$ is compact. We easily obtain the following solvability result of the inverse problem (compare with Theorem 2.2  of \cite{Seo}).
\begin{prop}
Suppose $\psi\in \mce'(\mbr^3)$ and $\WF(\psi)\cap \Sigma = \emptyset$. Then there exists $\chi\in \mcd'(\mbr^3)$ to
\beq
P(\p)\chi = -\lap \psi,
\eeq
such that $\singsupp(\chi) \subset \singsupp(\psi)$ is compact.
\end{prop}
The reconstruction formula of $\chi$ is obtained in Theorem 2.2 \cite{Seo}.

Based on the two propositions, the streaking artifacts can be identified as the set $ \WF(\chi)\backslash \WF(\psi)$. This is because for any $(x, \xi) \in \WF(\chi)\backslash \WF(\psi)$, the whole bicharacteristics $\gamma_{x, \xi}$ is in $\WF(\chi)$ and this is not compactly supported. The projection of $\gamma_{x, \xi}$ to $\mbr^3$ is the straight line from $x$ in $\xi$ direction. This also agrees with the numerical results, see Figure 1 of \cite{Seo}. Since the set $\WF(\chi)\backslash \WF(\psi)$ is contained in $(\La_p'\circ \WF(\psi))\backslash \WF(\psi)$, we shall regard the latter as the set of streaking artifacts in the following analysis.

\section{Reduction of the streaking artifacts}
For $\psi \in \mce'(\mbr^3)$, we decompose the solution $\chi$ to \eqref{eqwave} to separate the streaking artifacts.  Since the streaking artifacts are caused by $\WF(\psi)\cap \Sigma$, we decompose $\psi$ as following. Let $f\in C_0^\infty(\mbr)$ be a cut-off function such that $f(t) = 1, |t|<1$ and $f(t) = 0, |t|> 2$. Also, we let $g\in C_0^\infty(\mbr^3)$ be a cut-off function such that $\supp(\chi)\subset \supp(g)$. For $\eps >0$, we define a pseudo-differential operator $B_\eps$ of order $0$ associated with the (full) symbol 
\beqq\label{defBep}
b_\eps(x, \xi) = g(x)f(p(\xi)/\eps).
\eeqq
Here we recall that a pseudo-differential operator $A$ of order $m$ is defined by an oscillatory integral 
\beq
Au(x) = \frac{1}{(2\pi)^3}\int_{\mbr^3} e^{i(x-x')\xi} a(x, \xi) u(x') dx'd\xi, \ \ a(x, \xi)  \in S^m(\mbr^3\times \mbr^3),
\eeq
where $S^m(\mbr^3\times \mbr^3)$ denotes the standard symbol class i.e.\  for $a \in S^m(\mbr^3\times \mbr^3)$ and for any compact set $K$ of $U$, we have
\beq
|\p_x^\alpha \p_\xi^\beta a(x, \xi)| \leq C_{K, \alpha, \beta} \langle \xi \rangle^{m-|\beta|}, \ \ C_{K, \alpha, \beta} > 0,
\eeq
see e.g.\ \cite[Section 18.1]{Ho3}. $a$ is called the (full) symbol of $A$ and the principal symbol is defined in $S^m(\mbr^3\times \mbr^3)/S^{m-1}(\mbr^3\times \mbr^3).$  We denote the space of pseudo-differential operators of order $m$ by $\Psi^m(\mbr^3)$.  Observe that in \eqref{defBep} $b_\eps$ vanishes for $p(\xi)> 2\eps$ and $1-b_\eps$ vanishes for $p(\xi)< \eps$. Therefore, $\WF((\Id - B_\eps) \psi) \cap \Sigma = \emptyset.$ We can write
\beqq\label{decom}
\chi = \chi_1 + \chi_2, \text{ where } \chi_1 = Q(\Id - B_\eps) (-\lap \psi), \ \ \chi_2 = QB_\eps (-\lap \psi).
\eeqq
By Proposition \ref{wf}, we see that $\WF(\chi_1)\subset \WF(\psi)$ hence the streaking artifacts are contained in $\chi_2$. One can   remove the artifacts by simply taking $\chi_1$ as the reconstruction. However, this would remove all the singularities of $\chi$ on $\WF(\psi)\cap \Sigma$ and this would result in a loss of information. See Figure \ref{fig:TKD_X1_X2}.

To improve the results, in the following, we shall assume that the singularities of $\psi$ have special structure i.e.\ the singularities are in the normal directions of some submanifolds. These are called conormal distributions and they appear often in applications. For example, the delta distribution $\delta$. We see that $\singsupp(\delta) = \{0\in \mbr^3\}$ and $\WF(\delta) = T_0^*\mbr^3\backslash 0$. Here $0$ represents the zero section of $T^*\mbr^3$. Another example is the characteristic function $\mathbf{1}_\Omega$ where $\Omega$ is a bounded domain with smooth boundary $\p \Omega$. We have $\singsupp(\mathbf{1}_\Omega) = \p\Omega$ and $\WF(\mathbf{1}_\Omega) = \{(x, \xi)\in T^*\mbr^3: x\in \p\Omega, \langle \xi, \theta\rangle = 0, \forall\theta \in T_x(\p \Omega)\}$. The conormal bundle $N^*K$ of a submanifold $K\subset \mbr^3$, defined as
\beq
N^*K = \{(x, \xi)\in T^*\mbr^3 \backslash 0: \langle \xi, \theta \rangle = 0, \theta \in T_x K\},
\eeq
is a conic Lagrangian submanifold of $T^*\mbr^3\backslash 0$. In particular, $\WF(\delta) = N^*\{0\}$ and $\WF(\mathbf{1}_\Omega) = N^*(\p\Omega)$. 

We recall the basics of Lagrangian and paired Lagrangian distributions. The details can be found in for example \cite{DUV, Ho3}. Let $\La$ be a smooth conic Lagrangian submanifold of $T^*\mbr^3\backslash 0$. Following the standard notation, we denote by $I^\mu(\La)$ the space of Lagrangian distributions of order $\mu$ associated with $\La$. In particular, for $U$ open in $\mbr^3$, let $\phi(x, \xi): U\times \mbr^N \rightarrow \mbr$ be a smooth non-degenerate phase function (homogeneous of degree $1$ in $\xi$) that locally parametrizes $\La$ i.e. 
\beq
\{(x, d_x\phi)\in T^*_U \mbr^3 \backslash 0: x\in U,\ \ d_\xi \phi = 0\} \subset \La.
\eeq
Then $u\in I^{\mu}(\La)$ can be locally written as a finite sum of oscillatory integrals
\beq
\int_{\mbr^N} e^{i\phi(x, \xi)} a(x, \xi) d\xi, \ \ a\in S^{\mu + \frac 34 - \frac{N}{2}}(U\times \mbr^N).
\eeq
For $u\in I^\mu(\La)$, we know that $\WF(u)\subset \La$. Also, for any $s< -\mu-\frac{3}{4}$, we have  $u\in H^s(\mbr^3)$. So the order $\mu$ indicates the regularity of $u$. For a submanifold $Y\subset \mbr^3$, we denote $I^{\mu}(Y) = I^\mu(N^*Y)$ and these are called conormal distributions to $Y$. 

For two Lagrangians $\La_0, \La_1 \subset T^*X\backslash 0$ intersecting cleanly at a codimension $k$ submanifold i.e. 
\beq
T_x\La_0\cap T_x\La_1  = T_x(\La_0\cap \La_1),\ \ \forall x\in \La_0\cap \La_1,
\eeq
the paired Lagrangian distribution associated with $(\La_0, \La_1)$ is denoted by $I^{p, l}(\La_0, \La_1)$. Locally, paired Lagrangian distributions can be written as an oscillatory integral with a symbol of product type, see e.g.\ \cite{DUV}. For $u\in I^{p, l}(\La_0, \La_1)$, we know that $\WF(u)\subset \La_0\cup \La_1$. In particular, microlocally away from the intersection $\La_0\cap \La_1$, we have that $u \in I^{p+l}(\La_0\backslash \La_1)$ and $u \in I^p(\La_1\backslash \La_0)$. For example, we know that the fundamental solution $Q\in I^{-\frac 32, -\ha}(N^*\diag, \La_p)$ so that $Q\in I^{-2}(N^*\diag\backslash \La_p)$ and $Q\in I^{-\frac 32}(\La_p\backslash N^*\diag)$.

We need the following result to reduce the streaking artifacts. 
\begin{prop}\label{lmpsu}
Let $Y$ be a submanifold of $\mbr^3$ with codimension $\geq 1$ and assume that $N^*Y$ intersect $\Sigma$ transversally and each null bicharacterstics intersects $N^*Y$ a finite number of times. Let $\La_Y$ denote the flow out of $N^*Y$ under the Hamiltonian flow \eqref{flowout}. Let $K\in I^{p, l}(N^*\diag, \La_p)$ and $R\in \Psi^s(\mbr^3)$  properly supported. For $f\in I^\mu(Y)\cap \mce'(\mbr^3)$, we have
\begin{enumerate}
\item $R f \in I^{s+\mu}(Y)$. 
\item $K \circ R f \in I^{p+s + \mu, l}(N^*Y, \La_Y)$ so $K\circ Rf \in  I^{p+s + \mu+ l}(N^*Y\backslash \La_Y)\cap  I^{p+s + \mu}(\La_Y\backslash N^*Y)$.
\item  If the principal symbol of $R$ vanishes on  $\La_Y$, we have 
\beq
R\circ K f \in I^{p+\mu + l+s}(N^*Y\backslash \La_Y)\cap I^{p+\mu+ s-1}(\La_Y\backslash N^*Y).
\eeq
\end{enumerate}
\end{prop}

We make several remarks. First, one can think of the singularity on $N^*Y$ as the true singularity in $Kf$ and the one on $\La_Y$ as streaking artifacts. For $s<0$, part (2) says that one can reduce the singularity on the two pieces simultaneously. For $s =1$, part (3) tells one can increase the singularities on $N^*Y$ while the one on $\La_Y$ stays the same. Second,  if $R$ is of order $1$, the conclusion $R\circ K f \in I^{p+\mu}(\La_Y\backslash N^*Y)$ in part (3) just follows from the equivalent definition of Lagrangian distributions, see e.g.\ \cite[Definition 25.1.1]{Ho4}.

\bpf[Proof of Prop.\ \ref{lmpsu}]
(1) Let $k \geq 1$ be the codimension of $Y$. Locally we can choose local coordinate $x = (\bar x, \hat x), \bar x\in \mbr^{3-k}, \hat x \in \mbr^k$ such that $Y = \{\hat x = 0\}.$ We let the dual variables $\xi = (\bar \xi, \hat \xi),  \bar \xi \in \mbr^{3-k}, \hat \xi \in \mbr^k$. Then $N^*Y = \{\hat x = 0, \bar \xi = 0\}$. We can write $f\in I^\mu(Y)$ as an oscillatory integral
\beq
f(x) = \int_{\mbr^{k}} e^{i \hat x\cdot \hat \xi} b(\bar x, \hat \xi) d\hat \xi, \ \ b\in S^{\mu + \frac 34 - \frac{k}{2}}(\mbr^3\times \mbr^k).
\eeq
By partition of unity, we can assume that $f$ is compactly supported. For $R\in \Psi^s(\mbr^3)$ properly supported, we can write
\beq
\begin{split}
Rf(x) &= \frac{1}{(2\pi)^3}\int_{\mbr^3} e^{i (x-y)\cdot \xi} a(x, \xi)f(y) d\xi dy, \ \ a\in S^{s}(\mbr^3\times \mbr^3\times \mbr^3\backslash 0).
\end{split}
\eeq
Now we have
\beq
\begin{split}
Rf(x) =\frac{1}{(2\pi)^3} \int_{\mbr^3} \int_{\mbr^3} \int_{\mbr^{k}} e^{i (x-y)\cdot \xi} a(x, \xi)  e^{i \hat y\cdot \hat \eta} b(\bar y, \hat \eta) d\hat \eta d\xi dy =\frac{1}{(2\pi)^3} \int_{\mbr^k} e^{i\hat x \cdot \hat \xi} c(x, \hat \xi) d\hat \xi, 
\end{split}
\eeq
where
\beq
\begin{split}
 c(x, \hat \xi) &= \int_{\mbr^3}\int_{\mbr^{3-k}}\int_{\mbr^{k}} e^{i (\bar x-\bar y)\cdot \bar \xi} a(x, \xi)  e^{i \hat y\cdot (\hat \eta-\hat \xi)} b(\bar y, \hat \eta) d\hat \eta d\bar \xi dy.
\end{split}
\eeq 
We let $\hat \xi = \la \hat \theta$ with $\la = |\hat \xi|$ and let $\hat\alpha = \la^{-1}(\hat\eta - \hat \xi)$ i.e.\ $\hat \eta = \la\hat\alpha + \hat\xi$. Also, we let $\bar \beta = \la^{-1} \bar \xi$. So $\xi = (\bar \xi, \hat \xi) = \la (\bar \beta, \hat\theta)$. After these changes of variables, we get
\beq
\begin{split}
 c(x, \hat \xi) &= \int_{\mbr^3}\int_{\mbr^{3-k}}\int_{\mbr^{k}} e^{i \la (\bar x-\bar y)\cdot \bar \beta} a(x, (\la \bar \beta, \hat \xi))  e^{i \la \hat y\cdot \hat \alpha} b(\bar y, \la\hat\alpha + \hat\xi) \la^{3}d\hat \alpha d\bar \beta dy\\
 & = (-1)^{k}\int_{\mbr^3}\int_{\mbr^{3-k}} \int_{\mbr^{k}} e^{-i \la [\bar y\cdot \bar \beta + \hat y\cdot \hat \alpha]} a(x, (\la \bar \beta, \hat \xi))   b(\bar y + \bar x, \la\hat\alpha + \hat\xi) \la^{3}d\hat \alpha d\bar \beta dy.
\end{split}
\eeq 
Now we can apply stationary phase  on variables $y, \hat \alpha, \bar \beta$ to obtain 
\beq
c(x, \hat\xi) = C a(x, (0, \hat\xi) ) b(\bar x, \hat \xi) + \cdots,
\eeq
which belongs to $S^{\mu + s + \frac 34 - \frac{k}{2}}(\mbr^3\times \mbr^k)$. Here $C$ is a constant and the terms in $\cdots$ belong to $S^{\mu + s + \frac 34 - \frac{k}{2}-1}(\mbr^3\times \mbr^k)$, see for example the proof of \cite[Theorem 3.4]{GS}. Thus $Rf\in I^{\mu + s}(Y)$.

(2) We know that $Rf \in I^{\mu+s}(Y)$ from part (1). We can apply Prop.\ 2.1 of \cite{GrU1} to get $K\circ Rf \in I^{p + \mu+ s, l}(N^*Y, \La_Y)$. 

(3) We first apply Prop.\ 2.1 of \cite{GrU1} to conclude that $Kf\in I^{p+\mu, l}(N^*Y, \La_Y)$. Then we know that $Kf\in I^{p+\mu +l}(N^*Y\backslash \La_Y)$ and $Kf\in I^{p+\mu}(\La_Y\backslash N^*Y)$ as Lagrangian distributions. For $R\in \Psi^s(\mbr^3)$, we know (see e.g.\ Lemma 7.2 of \cite{GS}) that $\WF(Ru) \subset  \WF(u)$. Now we can apply part (1) to conclude that $R\circ Kf\in I^{p+\mu + s + l}(N^*Y \backslash \La_Y)$. If the  principal symbol of $R$ vanishes on $\La_Y$, we examine the proof of part (1) that the order of $R\circ Kf$ is $p+\mu-1+s$. This finishes the proof.
\epf

Now let's consider using part (2) of Lemma \ref{lmpsu} to reduce streaking artifacts. We take $R\in \Psi^{-s}(\mbr^3)$ with $s>0$, for example, 
\beq
Ru(x) = \frac{1}{(2\pi)^3}\int_{\mbr^3} e^{ix \cdot \xi} (1- f(|\xi|/\eps))|\xi|^{-s} \hat u(\xi) d\xi
\eeq
with $f$ the cut-off function defined at the beginning of this section. Then we replace $\chi_2$ by 
\beq
\tilde \chi_2 = Q\circ R\circ B_\eps (-\lap \psi),
\eeq
and consider an approximation $\tilde{\chi} = \chi_1 + \tilde{\chi}_2$ of $\chi$. Under the assumption of Lemma \ref{lmpsu}, if $\psi \in I^{\mu}(N^*Y),$ we know that $\tilde \chi_2 \in I^{-\frac 32 - s + 2+\mu, -\ha}(N^*Y, \La_Y)$. In particular, $\tilde \chi$ is a Lagrangian distribution on $N^*Y$ of order $\mu-s$ and a Lagrangian distribution of order $\mu+\ha -s$ on $\La_Y$. Recall that the streaking artifacts are contained in $\La_Y$, we conclude that by applying $R$, the singularities of streaking artifacts are reduced. However, this will also reduce the singularities on $N^*Y.$ 

Next we consider using part (3) of Lemma \ref{lmpsu} to relatively reduce the streaking artifacts. Let $R\in \Psi^s(\mbr^3), s >0$ with symbol vanishing on $\La_Y$. For example, we can take $R = P(\p)$ and replace $\chi_2$ by
\beq
 \hat \chi_2(x) = P(\p) \circ Q \circ B_\eps (-\lap \psi) = \frac{1}{(2\pi)^3}\int_{\mbr^3} e^{ix \xi} b_\eps(x, \xi)  |\xi|^2\hat \psi(\xi)d\xi.
\eeq
The key point is that the order of $\hat \chi_2 \in I^{-2+ \mu+s}(N^*Y\backslash \La_Y)$ is increased by $s=2$ $(s >0)$, while the singularities $\hat \chi_2 \in I^{-\frac 32 + \mu + s-1}(\La_Y\backslash N^*Y)$ is only increased by $s-1=1$. In other words, although the streaking artifacts are not reduced, the desired singularities in $\chi$ are enhanced. So the streaking artifacts are relatively reduced. 
As another example, we decompose the symbol 
\beq
p(\xi) = -\frac{1}{3} (\sqrt 2 \xi_3 - \sqrt{\xi_1^2 + \xi_2^2})(\sqrt 2 \xi_3 + \sqrt{\xi_1^2 + \xi_2^2}).
\eeq
Let $h(t)$ be a smooth cut-off function with $h(t) = 0, t< 0$. We let 
\beqq\label{op_T}
Ru(x) = \frac{1}{(2\pi)^3}\int_{\mbr^3}\int_{\mbr^3} e^{i(x-y)\xi} [h(\xi_3)(\sqrt 2 \xi_3 - \sqrt{\xi_1^2 + \xi_2^2}) + h(-\xi_3)(\sqrt 2 \xi_3 + \sqrt{\xi_1^2 + \xi_2^2})] u(y) dyd\xi.
\eeqq
In particular, $R\in \Psi^{1}(\mbr^3)$ and the symbol of $R$ vanish on $\La_Y$.

Finally, we can combine the above two approaches. Let $\psi\in I^\mu(Y)$. For $m$ a positive integer, we let $R\in \Psi^{-m}(\mbr^3)$ and $T\in \Psi^{1}(\mbr^3)$ so that the symbol of $T$ vanishes on $\La_Y$. 
Now we consider 
\beq
\chi_2^\# =  T^m \circ Q \circ R\circ B_\eps (-\lap \psi) 
\eeq
as a substitute for $\chi_2$. Then we have $\chi_2^\#\in I^{\mu}(N^*Y\backslash \La_Y)$ while $\chi_2^\#\in I^{\mu-m+\ha}(\La_Y\backslash N^*Y)$. So the streaking artifacts are reduced while the singularities on $N^*Y$ remain of the same strength.

\section{Numerical experiments}
\begin{figure}
\centering
\includegraphics[width=0.8\textwidth]{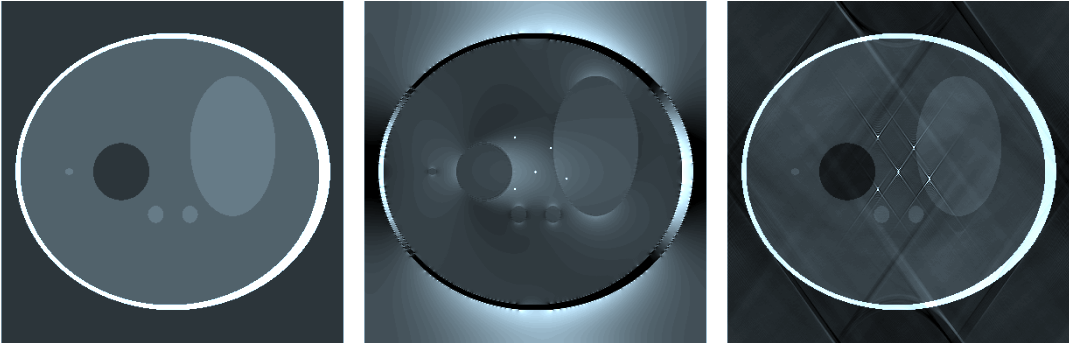}
\vspace{-.7em}
\fnote{(a)\hspace{11.5em}(b)\hspace{11em}(c)}
\caption{(a) Reference magnetic susceptibility distribution $\chi$. (b) Simulated data $\psi$ perturbed with point singularities. (c) Reconstructed $\chi$ using formula 2.18 in \cite{Seo}. The zero cone $\Sigma$ produces the propagation of singularities in the directions lying in the singular support of the fundamental solution $g(x)$ in \eqref{fund_sol}.}
\label{fig:phantom_data_seo}
\end{figure}

We carried out some numerical experiments using a grid of 392x392x392
with the purpose of illustrating the above theoretical analysis on the reduction of
streaking artifacts. The synthetic data was generated from the three dimensional
magnetic susceptibility distribution given by the
Shape-Loggan phantom following \eqref{qsm2}. As in \cite{Seo} we
perturbed the data by adding some point singularities in order to study their
propagation and aiming to reduce the streaking artifacts
caused by them (see figure \ref{fig:phantom_data_seo}). In what follows all the images correspond to sagittal views of three dimensional functions at $y=0$, and all the images corresponding to the magnetic susceptibility distribution are displayed using a window level of [-0.3, 1]. For the data $\psi$ we used the window level [-0.1,0.25].\\

One widely used reconstruction method in quantitative susceptibility mapping is the TKD method (see \cite{Seo} and references therein), which aims to recover $\chi$ by direct computations following the formula
$$\hat{\chi}_{\hbar}(\xi) = \left\{\begin{array}{ll} \frac{\hat{\psi}(\xi)}{D(\xi)}&\text{if } |D(\xi)|\geq\hbar,\\ \text{sign}(D(\xi))\frac{\hat{\psi}(\xi)}{\hbar}&\text{if } |D(\xi)|<\hbar.\end{array}\right.$$

For our experiments we implemmented a modification of the TKD method where the only difference is that we considered a smooth cut-off function to smoothly divide the Fourier Transform of the
reconstructed image into two pieces, the one supported away from
the characteristic set $\Sigma$ and the other one supported in a neighborhood of
it. Namely, the reconstructed image is of the form $\chi^{\hbar} =
\chi_1^{\hbar}+\chi_2^{\hbar}$ with
\beqq\label{TKD}
\chi_1^{\hbar} = Q\circ (\text{Id}-B_{\hbar})(-\Delta \psi),\quad
  \chi_2^{\hbar} = P_\hbar\circ Q\circ R\circ B_{\hbar}(-\Delta \psi),
\eeqq
where $Q$ and $B_{\hbar}$ are as above respectively in \eqref{defQ} and
\eqref{defBep} (with $g=1$); $P_\hbar\in\Psi^{2}(\mbr^3)$ given by the symbol $\sigma(P_\hbar)= \hbar^{-1}\text{sign}(p(\xi))p(\xi)$ which vanishes in $\Lambda_Y$, and $R\in\Psi^{-2}(\mbr^3)$ given by $\sigma(R) = |\xi|^{-2}$. In the Fourier domain  the
previous translates into the simpler formulas
\beq
\hat{\chi}^{\hbar}_1(\xi) = (1-b_\hbar(\xi))
\displaystyle\frac{\hat{\psi}(\xi)}{D(\xi)},\quad \hat{\chi}^{\hbar}_2(\xi) = b_\hbar(\xi)\text{sign}(p(\xi))
\displaystyle\frac{\hat{\psi}(\xi)}{\hbar}.
\eeq

\begin{figure}
\centering
\includegraphics[width=0.8\textwidth]{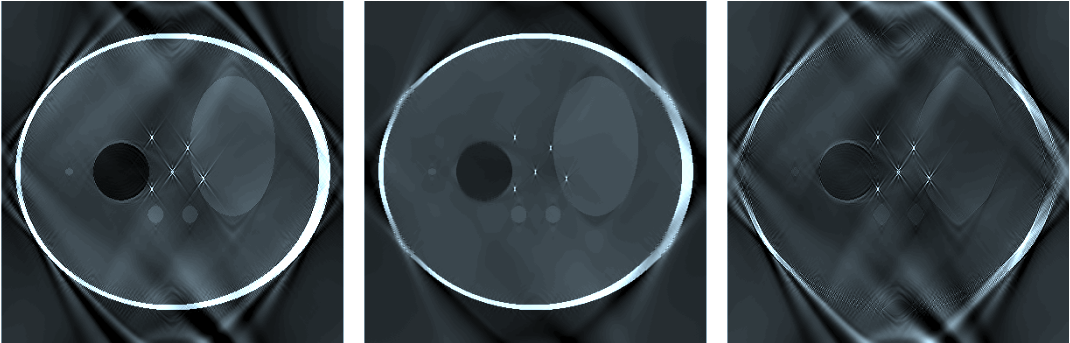}
\vspace{-.7em}
\fnote{(a)\hspace{11.5em}(b)\hspace{11em}(c)}
\caption{(a) Reconstructed susceptibility using the smooth version of the TKD method with $\hbar = 0.04$. (b) Reconstructed image using only frequencies away from the zero cone, i.e. $\chi^{\hbar}_1$. (c) Image obtained only considering frequencies near the zero cone, i.e. $\chi^{\hbar}_2$.}
\label{fig:TKD_X1_X2}
\end{figure}

\begin{figure}
\centering
\includegraphics[width=0.8\textwidth]{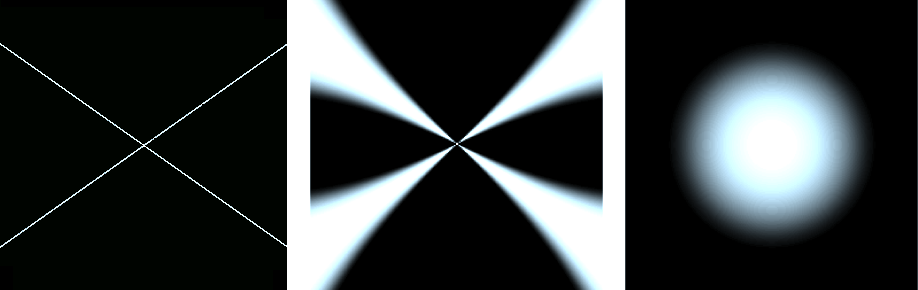}
\vspace{-.7em}
\fnote{(a)\hspace{11.5em}(b)\hspace{11em}(c)}
\caption{(a) Zero cone $\Sigma$. (b) Symbol of the cut-off pseudo-differential operator $B_{\hbar}$. (c) Symbol of the second cut-off operator $C_{M,\epsilon}$.}
\label{fig:Sigma_b_c}
\end{figure}

According to the theoretical analysis of the previous sections, what the TKD procedure does is first regularize all the singularities in directions near the zero cone by an order 2 and then enhance them, applying the operator $P_\hbar$, by and order of $2$ in $N^*Y\backslash \Lambda_Y$ and $1$ in $\Lambda_Y\backslash N^*Y$. Consequently, if the data belongs to $I^\mu(N^*Y)$ the streaking artifacts have order $\mu - 1/2$ and the other singularities have order $\mu$. As one can see in figure \ref{fig:TKD_X1_X2}(a), the streaking artifacts are slightly attenuated in contrast to figure \ref{fig:phantom_data_seo}(c) which was obtained by applying formula 2.18 in \cite{Seo}. The artifacts are still fairly visible though making necessary to use further methods to diminish them. As mentioned above, since the streaking artifacts are caused by the zero cone $\Sigma$ (figure \ref{fig:Sigma_b_c}(a)), removing the part of the image associated with such frequencies, this is only considering $\chi^\hbar_1$, implies lost of information as one can see in figure \ref{fig:TKD_X1_X2}(b) where some of the edges were smoothed out. The singularities appearing in $\chi^{\hbar}_1$ and $\chi^{\hbar}_2$ are sensitive to the shape of the symbol $b_\hbar$. Indeed, if we narrow the support of the symbol it causes the presence of streaking artifacts in $\chi^{\hbar}_1$. A sagittal view of the function $b_\hbar$ considered in our experiments is given in figure \ref{fig:Sigma_b_c}(b).

\begin{figure}
\centering
\includegraphics[width=0.8\textwidth]{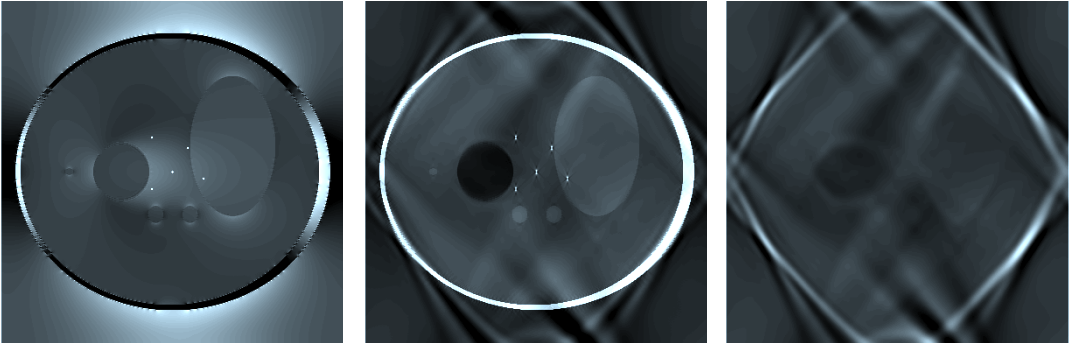}
\vspace{-.7em}
\fnote{(a)\hspace{11.5em}(b)\hspace{11em}(c)}
\caption{(a) Simulated data $\psi$. (b) Reconstructed image $\chi^{\hbar}$ applying the operator $R$ in \eqref{reg_op} with $s=2$ and $\hbar = 0.04$. (c) $\chi^\hbar_2$. The reduction of the streaking artifacts is clear and there is less lost of information in contrast to the case of just considering $\chi^\hbar_1$ (see figure \ref{fig:TKD_X1_X2}(b) ). In $(c)$ it can be seen that most of the smooth features contained in $\chi^\hbar_2$ are still there since $R$ only acts on $\chi^\hbar_{2,1}$.}
\label{fig:reg_ex}
\end{figure}

We can further reduce the order of the singularities as shown in figure \ref{fig:reg_ex} by increasing the order of the operator $R$. More precisely, we consider 
\beqq\label{reg_op}
R\in \Psi^{-s}(\mbr^3) \text{ with } s>2 \text{ and symbol }r(\xi) = K|\xi|^{-s},\; K>0.
\eeqq
Moreover, for some smooth function $f\in C^\infty_0(\mbr)$ such that $f(t) = 1$, $|t|<1$ and $f(t)=0$, $|t|>M$ for some $M>0$, we compute $\chi^\hbar_2 = \chi^\hbar_{2,1} + \chi^\hbar_{2,2}$, where
\beqq\label{R}
\chi^\hbar_{2,1} = P_\hbar\circ Q\circ R \circ (I-C)\circ B_\hbar(-\Delta\psi),\quad \chi^\hbar_{2,2} = P_\hbar \circ Q\circ C\circ B_\hbar(-\Delta\psi)
\eeqq
with $C=C_{M,\epsilon}\in \Psi^0(\mbr^3)$ of symbol $c_{M,\epsilon}(\xi) = f(|\xi|/\epsilon)$ for some $\epsilon>0$ and $M>0$ (see figure \ref{fig:Sigma_b_c}(c)). The motivation behind the division of $\chi^\hbar_2$ into two parts is that we would like to keep the smooth attributes of the image which are contained in $\chi^\hbar_{2,2}$, as well as reduce the artifacts included in $\chi^\hbar_{2,1}$. The images in figure \ref{fig:reg_ex} were obtained with $s=4$, therefore the streaking artifacts are reduced to order $\mu-1/2 - 2$ while the other singularities to order $\mu-2$.\\



\begin{figure}
\centering
\includegraphics[width=0.8\textwidth]{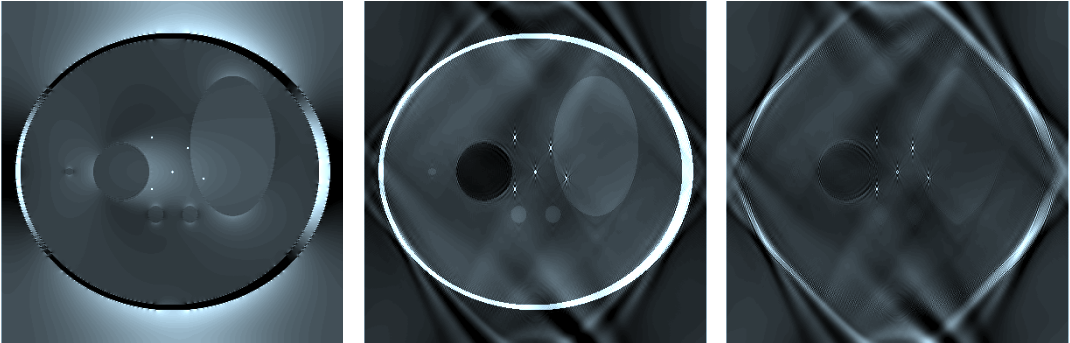}
\vspace{-.7em}
\fnote{(a)\hspace{11.5em}(b)\hspace{11em}(c)}
\caption{(a) Simulated data $\psi$. (b) Reconstructed susceptibility ${\chi}^{\hbar}$ following \eqref{T_op} with $m=s=2$ and $\hbar = 0.04$. (c) $\chi^{\hbar}_2$. The streaking artifacts were smoothed out but not the singularities related to edges of the phantom.}
\label{fig:T_ex}
\end{figure}

\begin{figure}
\centering
\includegraphics[width=0.8\textwidth]{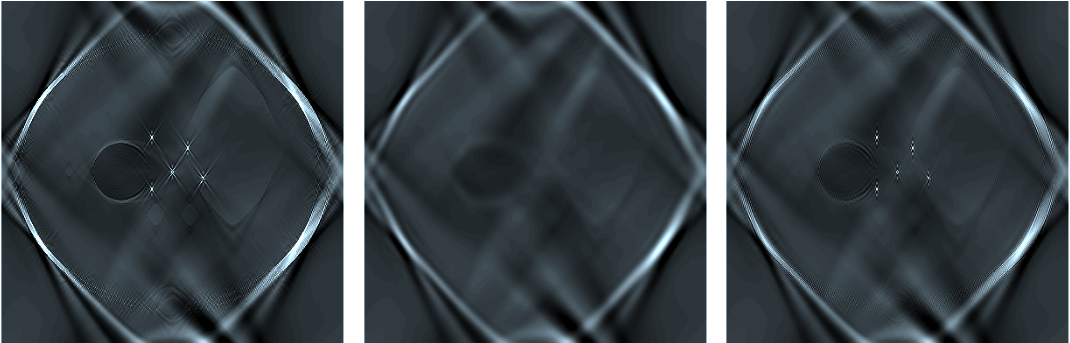}
\vspace{-.7em}
\fnote{(a)\hspace{11.5em}(b)\hspace{11em}(c)}
\caption{Comparison of $\chi_{2}^{\hbar}$ in the previous three cases: (a) smooth TKD, this is following \eqref{TKD}; (b) applying the operator $R$ as in \eqref{R} with $s=2$; (c) applying the operators $T$ and $R$ as in \eqref{T_op} with $m=2$ and $s=4$. Recall that the reconstructed susceptibility is given by $\chi^\hbar = \chi^\hbar_1 + \chi^\hbar_2$ with the last term containing the streaking artifacts.}
\label{fig:X2_comp}
\end{figure}

Finally, by considering the operator $T\in \Psi^1(\mbr^3)$ defined in \eqref{op_T} and $R\in\Psi^{-s}(\mbr^3)$ as in \eqref{reg_op}, the reconstructed susceptibility $\chi^{\hbar} = \chi^\hbar_1 + \chi^{\hbar}_2$ in figure \ref{fig:T_ex} is computed by doing $\chi^{\hbar}_2 = \chi^{\hbar}_{2,1} + \chi^{\hbar}_{2,2}$, where $m +2 = s > 0$ and
\beqq\label{T_op}
\chi^{\hbar}_{2,1} = (T^m\circ P_\hbar)\circ Q\circ R\circ (I-C)\circ B_\hbar(-\Delta\psi),\quad \chi^{\hbar}_{2,2} = P_\hbar \circ Q\circ C\circ B_\hbar(-\Delta\psi).
\eeqq
The streaking artifacts in $\chi^{\hbar}_{2,1}$ are reduced by an order of $-m-1/2$ while the rest of its singularities (the ones in $N^*Y\backslash  \Lambda_Y$) remain in the same order $\mu$ when $\psi\in I^\mu(N^*Y)$. By implementing this procedure the streaking artifacts are further reduced, in comparison with applying just TKD, and there is no attenuation of the rest of the singularities as for instance can be noticed in figure \ref{fig:X2_comp}, which shows the part of the reconstructed susceptibilities that contains the streaking artifacts, this is $\chi^\hbar_2$, for the previous three cases.

\end{document}